\documentclass[12pt]{amsart}
\usepackage[centertags]{amsmath}
\usepackage{amsfonts,amssymb}
\usepackage{graphicx}
\usepackage{enumerate}
\usepackage[latin1]{inputenc}

  \newtheorem{theorem}{Theorem}
  \newtheorem{corollary}{Corollary}
  \newtheorem{proposition}{Proposition}
  \newtheorem{lemma}{Lemma}%
  \theoremstyle{remark}

\begin{document}

\title{Coprime matchings}
\date{\today}

\author{Carl Pomerance}
\address{Mathematics Department, Dartmouth College, Hanover, NH 03784}
\email{carl.pomerance@dartmouth.edu}

\begin{abstract}
We prove that there is a matching between 2 intervals of positive integers of the same
even length, with corresponding pairs coprime, provided the intervals are in
$[n]$ and their lengths are $>c(\log n)^2$, for a positive constant $c$.
This improves on a recent result of Bohman and Peng.  As in their paper,
the result has an application to the lonely runner conjecture.
\end{abstract}



\maketitle

\section{Introduction}

Suppose one has two intervals $I$ and $J$ of positive integers, and both have
the same length.  Among the many possible bijections from $I$ to $J$,
it is interesting to wonder if at least one of them is a coprime matching.  That is,
in the bijection, corresponding numbers should be relatively prime.  It is easy to
see that sometimes this is not possible.  For example, if $1\notin I$ and $J$
contains a number $j$ divisible by the product of the members of $I$, then that
number $j$ cannot correspond to any member of $I$.  Another easy way to
block a coprime matching is if both intervals contain a strict majority of even
numbers, as would be the case with $\{2,3,4\}$ and $\{8,9,10\}$, for example.

An easy way to bar the majority evens case is to insist that the common length be
even, say $2m$.  And an easy way to bar an element from one interval having
a common nontrivial factor with each number of the other interval is to insist
that the numbers are not too large in comparison to the length.

Here are some existing results about coprime matchings of intervals.
If one of the intervals is $[n]=\{1,2,\dots,n\}$ and the other an arbitrary interval of length $n$, then there is a coprime matching, see \cite{PS}.  Here we allow
odd lengths and there's no prohibition of one of the intervals involving numbers
much larger than the length, but these are compensated for by having the
number 1 in one of the intervals.

Very recently, Bohman and Peng \cite{BP} 
showed that if $I,J$ are contiguous intervals of any length $k\ge 4$
such that $k\in I$,
then there is a coprime matching, so verifying a conjecture of
Larsen, Lehmann, Park, and Robertson \cite{LLPR}.
They also proved the following result.

\medskip
\noindent{\bf Theorem A}.  {\it 
There is a positive constant
${\mathcal C}$ such that if $n$ is sufficiently large, 
$m>\exp({\mathcal C}(\log\log n)^2)$, and
$I,J\subset[n]$
are intervals of length $2m$,  then there is a coprime matching
of $I$ and $J$.}
\medskip

In this note we improve on Theorem A.
\begin{theorem}
\label{thm:main}
There is a positive constant $c$ such that if $n$ is sufficiently large, 
$m>c(\log n)^2$, and $I,J\subset[n]$ are intervals of length $2m$,
then there is a coprime matching of $I$ and $J$.
\end{theorem}

Note that the expression $c(\log n)^2$ cannot be improved to $\log n$.
This follows from \cite[Theorem I]{E}, where it is shown that there are
infinitely many integers $n$ for which there is an interval of $\ge\log n$
integers in $[n]$ each with a nontrivial common divisor with $n$.
This can be improved a little using better estimates for the Jacobsthal
function, for example, \cite{F}.

Adding to the interest of the Bohman--Peng paper is an application of Theorem A
to the lonely runner conjecture.  
This conjecture asserts that
if $v_1,v_2,\dots, v_n$ are distinct positive integers, then there is a real number
$t$ such that no $v_it$ is within distance $1/(n+1)$ of an integer.  This
has been shown by Tao \cite{T} when the $v_i$'s are at most $1.2n$.  In \cite{BP},
the lonely runner conjecture is shown when the $v_i$'s are at most $2n-\epsilon(n)$, where
$\epsilon(n)$ is of the shape $\exp({\mathcal C}'(\log\log n)^2)$, for a positive
constant ${\mathcal C}'$.  
Theorem \ref{thm:main} (more precisely, Proposition \ref{prop:main} below)
has the analogous application: by the same argument
as in \cite{BP},
if $v_1<v_2<\dots<v_n\le 2n-c'(\log n)^2$, the lonely runner conjecture holds.
It remains a challenge to
show it when $\{v_1,v_2,\dots,v_n\}\subset[2n]$,  much less the full conjecture.

\section{The set up}

We say two integers $s,t$ are 2-coprime if no odd prime divides both $s$ and $t$,
that is, $\gcd(s,t)$ is a power of 2.
As in the Bohman--Peng paper, the problem can be reduced to
the following.
\begin{proposition}
\label{prop:main}
There is a positive constant $c$ such that if $n$ is sufficiently large,
$m>c(\log n)^2$, and
 $I, J\subset[n]$ are arithmetic progressions of length $m$ with 
common difference $1$ or $2$, the following holds.
Whenever $S\subset I$ and $T\subset J$ are nonempty, 
with $|S|+|T|\ge m$, there is a $2$-coprime
pair $s,t$ with $s\in S$, $t\in T$. 
\end{proposition}
\begin{corollary}
\label{cor:hall}
There is a matching of $I$ and $J$ with corresponding numbers being $2$-coprime.
\end{corollary}
\begin{proof}
This follows immediately from Proposition \ref{prop:main} and Hall's matching
theorem. 
\end{proof}
As pointed out in \cite{BP}, to use Hall's theorem it suffices to consider the
case when $|S|+|T|>m$.  The weaker hypothesis $|S|+|T|\ge m$ in Proposition \ref{prop:main}
is useful when considering the application to the lonely runner problem.
\begin{corollary}
There is a positive constant $c$ such that if $n$ is sufficiently large,
$m>c(\log n)^2$, and $I,J\subset[n]$ are intervals, then there is a coprime
matching of $I$ and $J$ if either 
\begin{enumerate}
\item $I,J$ have length $2m$, or 
\item
$I,J$ have length $2m+1$ and the least elements of $I,J$ have opposite parity.
\end{enumerate}

\end{corollary}
\begin{proof}
First suppose that $I,J$ have length $2m$.  Let $I_0,J_0$ be the even elements of $I,J$,
respectively, and let $I_1,J_1$ be the odd elements.  We have
$|I_0|=|J_1|=m$ and $|I_1|=|J_0|=m$.  Applying Corollary \ref{cor:hall}, there
is a matching of $I_0$ and $J_1$ with corresponding elements being 2-coprime.
But as elements of $I_0$ are all even and elements of $J_1$ are all odd,
being 2-coprime implies being coprime.  So the matching is a coprime one,
and similarly for $I_1$ and $J_0$.  Thus, we have the result in the first case.
So suppose $I,J$ have length $2m+1$ with least elements of the opposite
parity.  Then again we have $|I_0|=|J_1|$ and $|I_1|=|J_0|$,
and the same proof works.
\end{proof}
In particular, Theorem \ref{thm:main} holds.  Probably there should be a coprime
matching in the length $2m+1$ case when $I,J$ both have odd least elements.
In the length $2m+1$ case when both $I,J$ start with even numbers, perhaps
there is a matching where every pair but one is coprime, and the offending
pair has gcd 2.

\subsection{Sketch of the proof}

The argument in \cite{BP} uses a result of Erd\H os \cite{E} 
on the Jacobsthal function
that implies that a long string of consecutive members of an arithmetic progression
of common difference 1 or 2 has at least one member 2-coprime to a
given integer $s$.  We use instead a sequel result 
 of Iwaniec \cite{I} that implies that each $s\in S$
is 2-coprime to many elements of $J$, in fact, so many elements that we are
done unless $T$ is small enough to miss all of them.  But $T$ small forces
$S$ to be somewhat large, namely at least of magnitude $m/(\log m)^2$.
At that point, an averaging argument comes into
play.  In particular, it is first shown that such a large set $S$ has
most members with a not-too-large
``$m$ part" (namely the largest squarefree divisor composed of odd primes 
up to $m$).  This, coupled with the first argument shows that
we may assume that $|S|$ is even larger, at
least of magnitude $m/(\log\log m)^2$.  A finer averaging argument now shows
that for most $s\in S$ the value of $\varphi(s)/s$ is mostly determined by
the primes $\le\log m$ dividing $s$.  A final averaging argument
shows that many elements of $S$ have $\varphi(s)/s$ not too small.
Returning to the first thought of
getting many members of $J$ to be 2-coprime to $s$, we no longer need the
Iwaniec result and can do a complete inclusion-exclusion, which allows us
to complete the proof.  This last step has some overlap with the approaches in
\cite{BP} and \cite{PS}.

\section{The proof of Proposition \ref{prop:main}}

Let $\varphi$ denote Euler's function, and let $\omega(n)$ denote the
number of different prime numbers that divide $n$.  For $n>2$ we have
$\omega(n)=O(\log n/\log\log n)$.  It is convenient to have a weaker, but
explicit inequality:  if $n>1$,
then $\omega(n)\le2\log n$ (since $\omega(n)$ is at most the base-2 logarithm
of $n$).  Further, for $\omega(n)>1$ and $n$ odd, $n/\varphi(n)\le3\log\omega(n)$.
This follows from considering those $n$ that are the product of the first $k\ge2$
odd primes and using (3.5) and (3.30) from \cite{RS}.

Assume the hypotheses of Proposition \ref{prop:main} hold.  Let $S\subset I$,
$T\subset J$ be nonempty subsets with $|S|+|T|\ge m$.
We may assume that $|S|+|T|=m$ and $|S|\le |T|$.  
For any integer $k>0$, let $k_m$ denote
the largest odd, squarefree divisor of $k$ supported on the primes $\le m$.
First note that if $S$ contains an element $s$ such that $s_m=2^aq^b$ where $q$ is
an odd prime $\le n$ and $a,b\ge0$, then the number of elements of $J$ that are
not 2-coprime to $s_m$ is $\le m/q+1$.   Since primes $p>m$ can divide
at most 1 element of $J$, the number of elements of $J$ that are not 2-coprime
to $s$ is at most $m/q+1+\omega(s)$.  Since $\omega(s)\le2\log m$, it follows that
for large $m$, $J$ contains $<m/2$ numbers not 2-coprime to $s$.  
Since $|T|\ge m/2$ it thus must have a number 2-coprime to $s$.

From now on, we assume that each element of $S$ is divisible by at least two
different odd prime numbers $\le m$.   Let $s\in S$.  From \cite{I}, 
 it follows that  there is a positive constant $c_1$,
 such that an interval of length
$c_1(s_m/\varphi(s_m))\omega(s_m)^2\log\omega(s_m)$ has $\ge \omega(s_m)^2$ integers
coprime to $s_m$.  If $J$ has common
difference 1, we apply this result directly to sub-intervals of $J$.  If $J$ has
common difference 2 and $J\subset2{\mathbb Z}$, we apply it to $\frac12J$,
while if $J\subset2{\mathbb Z}+1$, we apply it to $\frac12(J+M)$, where $M$
is the product of all of the odd primes $\le m$.  In all cases we thus have that
for each string of $c_1(s_m/\varphi(s_m))\omega(s_m)^2\log\omega(s_m)$
consecutive members of $J$, there are at least $\omega(s_m)^2$ integers
coprime to $s_m$.

Note that $s_m$ satisfies
$s_m/\varphi(s_m)\le3\log\omega(s_m)$ as remarked above.
 Thus, a string  of length
$3c_1(\omega(s_m)\log \omega(s_m))^2$ of consecutive members of $J$ has
at least $\omega(s_m)^2$ numbers coprime to $s_m$.  
Further, since $\omega(s_m)=O(\log s_m/\log\log s_m)$, we have
$\omega(s_m)\log\omega(s_m)\le c_2\log s_m\le c_2\log n$ for some absolute constant $c_2>0$.
We let $c$ in the theorem
be $3c_1c_2^2$, so that $J$ has at least 1  string  of length $3c_1(\omega(s_m)\log\omega(s_m))^2$ of consecutive members.  In particular,
breaking $J$ (or $\frac12J$ in the case that $J\subset 2{\mathbb Z}$, or
$\frac12(J+M)$ in the case $J\subset2{\mathbb Z}+1$)
 into consecutive strings of length $3c_1(\omega(s_m)\log \omega(s_m))^2$, 
we see that $J$ contains
 $> (1/6c_1)m/(\log\omega(s_m))^2$ integers coprime to $s_m$.  As above, $s$
 is divisible by at most $2\log n$ primes larger than $m$, and each of these
 primes divides at most one member of $J$.  So $J$ contains at least
 $(1/6c_1)m/(\log\omega(s_m))^2-2\log n$ numbers 2-coprime to $s$.
Since $m> c(\log n)^2$ and $\omega(s_m)<2\log n$, we have
 that there is a positive constant $c_3$ such that for each $s\in S$,
 \begin{equation}
 \label{eq:Jbound}
 \sum_{\substack{j\in J\\j\,2\hbox{-}{\rm coprime\, to}\,s}}1\ge \frac{c_3m}{(\log\omega(s_m))^2}.
 \end{equation} 
 Hence we may assume that $|T| \le m-c_3m/(\log\omega(s_m))^2$, so that
 \begin{equation}
 \label{eq:Sbound}
 |S|\ge c_3m/(\log\omega(s_m))^2.
  \end{equation}
  Now $\log\omega(s_m)\le\log(3\log n)$ and $m>c(\log n)^2$ so that
  $\log\log n=O(\log m)$.  Thus \eqref{eq:Sbound} implies that there is a
  positive constant $c_4$ with
  \begin{equation}
  \label{eq:Sbound2}
  |S|\ge c_4m/(\log m)^2.
  \end{equation}
 
 \begin{lemma}
 \label{lem:smlarge}
 The number of integers $i\in I$ with $i_m>\exp({(\log m)^4})$ is $O(m/(\log m)^3)$.
 \end{lemma}
 \begin{proof}
 We have 
 \[
 \sum_{i\in I}\log i_m=\sum_{2<p\le m}\log p\sum_{\substack{i\in I\\p\,\mid\, i}}1
 \le2m\sum_{p\le m}\frac{\log p}p=O(m\log m),
 \]
 by an inequality of Chebyshev.  Thus, the lemma follows.
 \end{proof}
 
Using \eqref{eq:Sbound2}, Lemma \ref{lem:smlarge}
implies that there is a member
$s$ of $S$ with $s_m\le e^{(\log m)^4}$.  Thus, by \eqref{eq:Sbound},
we now can assume that
\begin{equation}
\label{eq:lb}
|S|\ge c_5\frac{m}{(\log\log m)^2}.
\end{equation}
  
 \begin{lemma}
 \label{lem:slogmlarge}
 The number of integers $i\in I$ with
 \[
 f(i):=\sum_{\substack{p\,\mid\, i\\\log m < p\le m}}\frac1p > \frac{(\log\log m)^4}{\log m}
 \]
 is at most $O(m/(\log\log m)^3)$.
 \end{lemma}
 \begin{proof}
  We have
 \begin{align*}
 \sum_{i\in I}f(i)&=\sum_{\log m<p\le m}\frac1p\sum_{\substack{i\in I\\p\,\mid\, i}}1
 \le\sum_{\log m <p\le m}\frac1p\left(\frac mp+1\right)\\
 &\le2m\sum_{\log m < p\le m}\frac1{p^2}=O\left(\frac m{\log m\log\log m}\right).
 \end{align*}
 The lemma follows.
 \end{proof}
 
 By \eqref{eq:lb}, Lemmas \ref{lem:smlarge} and
  \ref{lem:slogmlarge} imply
 that we may assume that there are at least $(1-2/(c_5\log\log m))|S|$
 members $s$ of $S$ with
 $s_m\le e^{(\log m)^4}$ and $f(s)\le (\log\log m)^4/\log m$.   For a positive
 integer $k$, let
 \[
 k_0=k_{\log m}=\prod_{\substack{p\,\mid\, k\\2<p\le\log m}}p.
 \]
 For any $s$ as above, the number of members of $J$ coprime to $s_0$ 
 is 
 \[
 \sum_{d\,\mid\, s_0}\mu(d)\sum_{\substack{j\in J\\d\,\mid\, j}}1\ge\sum_{d\,\mid\, s_0}
\left( \mu(d)\frac md-1\right)=\frac{\varphi(s_0)}{s_0}m-2^{\omega(s_0)}.
\]
Now $\omega(s_0)<\pi(\log m)=o(\log m)$, so that $2^{\omega(s_0)}=m^{o(1)}$,
and on the other hand, $m\varphi(s_0)/s_0=\Omega(m/\log\log m)$.
Thus, $J$ contains at least
 $0.99m\varphi(s_0)/s_0$ integers coprime to $s_0$.  The number of members
 of $J$ coprime to $s_0$ but not coprime to $s_m$ is at most 
 \[
 \sum_{\substack{p\,\mid\, s\\\log m<p\le m}}\left(\frac mp+1\right)\le2mf(s).
 \]
Since $f(s)\le(\log\log m)^4/\log m$, $J$ contains at least 
 $0.98m\varphi(s_0)/s_0$ integers coprime to $s_m$, and thus, as above, at least 
 $0.97m\varphi(s_0)/s_0$ integers 2-coprime to $s$.  And, as we have seen,
 this holds for at least $(1-2/(c_5\log\log m))|S|$ members $s$ of $S$.
 
 We use the next result (cf.\ \cite[Prop. 3]{PS})
 to show that it is unusual for an element $i\in I$ to
 have $i_0/\varphi(i_0)$ large.  
 \begin{lemma}
\label{lem:phi}
For $m$ sufficiently large, we have
\[
\sum_{i\in I}\left(\frac{i_0}{\varphi(i_0)}-1\right)<\frac3{10}m.
\]
\end{lemma}
\begin{proof}
We have 
\[
\frac k{\varphi(k)}=\sum_{d\,\mid\, k}\frac{\mu(d)^2}{\varphi(d)}.
\]
Let $P$ denote the product of the odd primes $p\le\log m$.  Thus,
\begin{align}
\label{eq:2p}
\sum_{i\in I}\frac {i_0}{\varphi(i_0)}&=\sum_{i\in I}
\sum_{d\,\mid\, i_0}\frac{\mu(d)^2}{\varphi(d)}
=\sum_{d\,\mid\, P}\frac1{\varphi(d)}\sum_{\substack{i\in I\\d\,\mid\, i}}1\notag\\
&\le\sum_{d\,\mid\, P}\frac{m/d+1}{\varphi(d)}=m\sum_{d\,\mid\, P}\frac 1{d\varphi(d)}+\sum_{d\,\mid\, P}\frac1{\varphi(d)}\notag\\
&=m\prod_{p\,\mid\, P}\left(1+\frac1{p(p-1)}\right)+\prod_{p\,\mid\, P}\,\frac p{p-1}.
\end{align}
Note that
\[
\prod_{p}\left(1+\frac1{p(p-1)}\right)=\frac{\zeta(2)\zeta(3)}{\zeta(6)}<1.944,
\]
where $\zeta$ is the Riemann zeta-function.  The first product in \eqref{eq:2p}
is missing the prime 2, so it is
$<(2/3)1.944=1.296$.  Extending the second product 
in \eqref{eq:2p} over all integers in
$[2, \log m]$, we see that it is $\le\log m$.  Thus, 
\[
\sum_{i\in I}\frac {i_0}{\varphi(i_0)}<1.296m+\log m,
\]
and the lemma follows for all sufficiently large $m$.
\end{proof}
\begin{corollary}
\label{cor:phi}
For all sufficiently large integers $m$ and for any
real number $t>1$, the number of $i\in I$ with $i_0/\varphi(i_0)>t$ is
at most $0.3m/(t-1)$.
\end{corollary}
\begin{proof}
Let $N$ denote the number of integers $i$ in question.  Lemma \ref{lem:phi} implies that
$N(t-1)<3m/10$.
\end{proof}

Let $r=m/|S|$, so that $r\ge2$.  We apply Corollary \ref{cor:phi} with $t=0.9r$,
and we deduce that the number of $i\in I$ with $i_0/\varphi(i_0)>0.9r$ is
$\le0.3m/(0.9r-1)\le3m/(8r)=\frac38|S|$.  Thus, more than $\frac58$ of the members $s$ of
$S$ have $s_0/\varphi(s_0)\le0.9r$.  Further, we have seen above that at least
$(1-2/(c_5\log\log m))|S|$ members $s$ of $S$ are 2-coprime to
 at least $0.97m\varphi(s_0)/s_0$ members of $J$.  Thus, at least 
 $(5/8-2/(c_5\log\log m))|S|$ members $s$ of $S$ have this property and also 
 $s_0/\varphi(s_0)\le 0.9r$.  Let $s$ be one such element.  There are at least
 $0.97m/(0.9r)$ elements $j$ of $J$ that are 2-coprime to $s$.  Now
 $0.97m/(0.9r)>1.07m/r=1.07|S|$, and $|T|=m-|S|$.  Thus some of these values of
 $j$ must be in $T$, completing the proof of Proposition \ref{prop:main}.


\end{document}